\documentclass[11pt]{article}
\usepackage{graphicx} % Required for inserting images
\usepackage[utf8]{inputenc}
\usepackage{titlesec}
\titleformat{\section}
{\normalfont\scshape}
{\thesection.}{0.5em}{\centering}
\titleformat{\subsection}
{\small\scshape}
{\thesubsection.}{0.5em}{\centering}
\titleformat{\subsubsection}
{\small\scshape}
{\thesubsubsection.}{0.5em}{\centering}
\usepackage{tikz-cd}
\usepackage{tikz}
\usepackage[left=3.5cm,right=3.5cm,top=3.5cm,bottom=3.5cm]{geometry}
\usepackage{comment}
\usepackage{amsfonts}
\usepackage{indentfirst}
\renewcommand{\abstract}{\small{\section*{\abstractname}}}
\usepackage{amsmath}
\usepackage{amssymb}
\usepackage{amsthm}
\usepackage{mathrsfs}
\newtheorem{theorem}{Theorem}
\newtheorem{lemma}{Lemma}

\newtheorem{corollary}{Corollary}

\newtheorem*{theorem*}{Theorem}
\theoremstyle{remark}
\newtheorem{remark}{Remark}
\usepackage{stmaryrd}

\usepackage[
backend=biber,
style=alphabetic,
sorting=ynt
]{biblatex}
\title{\Large{\textbf{Fatou limits of stochastic integrals}}}
\author{\normalsize{\textsc{Vasily Melnikov}}}
\date{\normalsize{\textsc{March 2025}}}

\begin{document}

\maketitle
\begin{abstract}
     The convergence of stochastic integrals is essential to stochastic analysis, especially in applications to mathematical finance, where they model the gains associated with a self-financing strategy. However, Fatou convergence of $(X^{n})_{n=1}^{\infty}$—a notion introduced for its amenability to compactness principles—implies little about the sequence of Itô integrals $\left(\int_{0}^{\cdot}YdX^{n}\right)_{n=1}^{\infty}$ for a fixed integrand $Y$. Under a boundedness condition, we find convex combinations $(\widetilde{X}^{n})_{n=1}^{\infty}$ of $(X^{n})_{n=1}^{\infty}$ with Fatou limit $\widetilde{X}$, such that $\left(\int_{0}^{\cdot}Yd\widetilde{X}^{n}\right)_{n=1}^{\infty}$ converges in a Fatou-like sense to $\int_{0}^{\cdot}Yd\widetilde{X}$ for all continuous semimartingales $Y$. The result is sharp, in the sense that continuity of $Y$ cannot be relaxed to being the left limits process of a semimartingale.
\end{abstract}
\section{Introduction}\label{sec:intro}
In their work on optional decompositions under constraints, Föllmer and Kramkov \cite{optdecompconstrain} defined the notion of a Fatou limit: given a dense subset $\mathscr{T}\subset[0,\infty)$, a sequence $(X^{n})_{n=1}^{\infty}$ of processes Fatou converge to $X$ on $\mathscr{T}$ if
\begin{equation*}
    X=\limsup_{s\downarrow\cdot,s\in\mathscr{T}}\limsup_{n\to\infty}X^{n}_{s}=\liminf_{s\downarrow\cdot,s\in\mathscr{T}}\liminf_{n\to\infty}X^{n}_{s}.
\end{equation*}
Originally introduced to prove hedgeability of certain contingent claims, Fatou limits have become a ubiquitous tool in mathematical finance—especially so in portfolio optimization, where they play a key role both in the frictionless case (see \cite{kramsch}), and in markets with transaction costs (see \cite{shadow}). Underlying its adoption is the relative weakness of Fatou convergence; it satisfies various compactness principles (see, for example, Lemma 5.2, \cite{optdecompconstrain}), making it easily applicable to a variety of situations.
\par
However, this relative weakness has a pernicious side; if $(X^{n})_{n=1}^{\infty}$ Fatou converges to $X$, it is not clear that the stochastic integrals $\left(\int_{0}^{\cdot}YdX^{n}\right)_{n=1}^{\infty}$ Fatou converge to $\int_{0}^{\cdot}YdX$ even for quite regular processes $Y$. Compactness principles for stochastic integrals, given their importance to mathematical finance (see \cite{convgemery}), have some precedent in the literature; however, the class of integrands deemed admissible is usually small, restricting to finite variation integrands (see \cite{martconvcomp}), or assuming the support of the integrand lies in some \textit{a priori} unknowable sequence of predictable sets (see \cite{del-sch} and \cite{melnikovemery}). We correct both of these issues, essentially proving that passing to Fatou limits for a sequence satisfying a boundedness condition preserves stochastic integration with respect to continuous semimartingale integrands—a wide class of processes.
\par
More precisely, our main result is the following. If a sequence $(X^{n})_{n=1}^{\infty}$ of semimartingales satisfies a boundedness condition, there exists $\widetilde{X}^{n}\in\mathrm{co}\{X^{m}:m\geq n\}$ and a semimartingale $\widetilde{X}$ such that $\widetilde{X}=\lim_{s\downarrow\cdot,s\in\mathbb{Q}_{+}}\lim_{n\to\infty}\widetilde{X}^{n}_{s}$ and
    \begin{equation*}
        \lim_{s\downarrow\cdot,s\in\mathbb{Q}_{+}}\mathbb{P}-\lim_{n\to\infty}\int_{0}^{s}Yd\widetilde{X}^{n}=\int_{0}^{\cdot}Yd\widetilde{X}
    \end{equation*}
    for each continuous semimartingale $Y$, where $\mathbb{P}-\lim_{n\to\infty}\int_{0}^{s}Yd\widetilde{X}^{n}$ denotes the limit in probability of $\left(\int_{0}^{s}Yd\widetilde{X}^{n}\right)_{n=1}^{\infty}$ for $s\in\mathbb{Q}_{+}$.
    \par
The notion of convergence above technically need not imply Fatou convergence; however, once one fixes $Y$, it is always possible to pass to a subsequence for which \textit{bona fide} Fatou convergence holds. Similar issues arose in previous works (see, for example, Proposition 4.1, \cite{martconvcomp}).
\par
The paper is structured as follows. In \S\ref{sec:notation}, we introduce our notation. In \S\ref{sec:main}, we state our main theorem, as well as some corollaries. In \S\ref{sec:proof}, we prove our main theorem. Finally, in \S\ref{sec:counter}, we provide counterexamples showing the sharpness of our results and techniques.
\section{Notation}\label{sec:notation}
Let $(\Omega,\mathscr{F},\mathbb{P})$ be a probability space. Suppose $(\xi_{n})_{n=1}^{\infty}$ is a sequence of random variables on $(\Omega,\mathscr{F},\mathbb{P})$. The limit in probability of $(\xi_{n})_{n=1}^{\infty}$, if it exists, will be denoted $\mathbb{P}-\lim_{n\to\infty}\xi_{n}$. The $\mathbb{P}$-a.s. limit of $(\xi_{n})_{n=1}^{\infty}$, if it exists, will be denoted $\lim_{n\to\infty}\xi_{n}$.
\par
Let $\mathbb{F}=\{\mathscr{F}_{t}:t\in[0,\infty)\}$ be a filtration of sub-$\sigma$-algebras of $\mathscr{F}$, satisfying the usual conditions, and such that $\mathscr{F}_{0}$ is the $\mathbb{P}$-augmentation of the trivial $\sigma$-algebra. Notions understood relative to a filtration (e.g. adaptedness, predictability) are understood to be relative to $\mathbb{F}$. If $X$ is a semimartingale, and $\xi$ is an $X$-integrable predictable process, the Itô integral of $\xi$ with respect to $X$ will be denoted by $\int_{0}^{\cdot}\xi dX$. If $X$ and $Y$ are semimartingales, the quadratic covariation of $X$ and $Y$ will be denoted by $[X,Y]$.
\par
Over the course of this paper, we will equip spaces of stochastic processes with various topologies. On the space of $\mathbb{F}$-adapted càdlàg processes modulo evanescence, we will consider the u.c.p. (shorthand for uniform convergence on compacts in probability) topology, defined by the translation-invariant metric
\begin{equation*}
    \mathbf{D}^{ucp}(X,Y)=\sum_{n=1}^{\infty}\frac{1}{2^{n}}\int_{\Omega}(X-Y)^{\ast}_{n}\wedge1d\mathbb{P}
\end{equation*}
where $X^{\ast}=\sup_{s\leq\cdot}\vert{X_{s}}\vert$ denotes the maximal function of $X$. On the space of semimartingales modulo evanescence, we will consider the Émery topology, defined by the translation-invariant metric
\begin{equation*}
   \mathbf{D}^{sm}(X,Y)=\sum_{n=1}^{\infty}\frac{1}{2^{n}}\sup_{\vert{\xi}\vert\leq1}\int_{\Omega}\left(\int_{0}^{\cdot}\xi d(X-Y)\right)^{\ast}_{n}\wedge1d\mathbb{P}. 
\end{equation*}
\par
If $f:[0,\infty)\longrightarrow\mathbb{R}$ is làdlàg, we will denote by $f_{-}$ the left limit of $f$ (with the convention that $f_{0-}=f_{0}$), and by $f_{+}$ the right limit of $f$. The left jump of $f$ is denoted by $\Delta f=f-f_{-}$. For later use, we will state Froda's theorem (in a slightly more general form, applying to finite variation functions of time).
\begin{lemma}\label{lem:froda}
    Suppose $f$ has locally finite variation. Then $f$ is làdlàg and
    \begin{equation*}
        \{t\in[0,\infty):f\textrm{ is discontinuous at }t\}=\{f\neq f_{-}\}\cup\{f\neq f_{+}\}.
    \end{equation*}
    Furthermore, the above sets are all at most countable.
\end{lemma}
\section{Main results}\label{sec:main}
The main theorem of this article, Theorem \ref{thm:main}, can be stated as follows.
\begin{theorem}\label{thm:main}
    Let $(X^{n})_{n=1}^{\infty}$ be a sequence of semimartingales such that, for each $t\geq0$, the set
    \begin{equation*}
        \mathrm{co}\left\{\left\vert{\int_{0}^{t}\xi dX^{n}}\right\vert:n\in\mathbb{N},\xi\textrm{ is predictable and }\vert{\xi}\vert\leq1\right\}
    \end{equation*}
    is bounded in probability, and $(X^{n}_{0})_{n=1}^{\infty}$ is bounded. Then there exists $\widetilde{X}^{n}\in\mathrm{co}\{X^{m}:m\geq n\}$ and a semimartingale $\widetilde{X}$ such that $\widetilde{X}=\lim_{s\downarrow\cdot,s\in\mathbb{Q}_{+}}\lim_{n\to\infty}\widetilde{X}^{n}_{s}$ and
    \begin{equation*}
        \lim_{s\downarrow\cdot,s\in\mathbb{Q}_{+}}\mathbb{P}-\lim_{n\to\infty}\int_{0}^{s}Yd\widetilde{X}^{n}=\int_{0}^{\cdot}Yd\widetilde{X}
    \end{equation*}
    for each continuous semimartingale $Y$.
\end{theorem}
The continuity condition on $Y$ situates Theorem \ref{thm:main} as a stochastic counterpart to the following well-known version of Prokhorov's theorem: if $(g_{n})_{n=1}^{\infty}$ is a sequence of finite variation functions on $[0,1]$ with $\sup_{n}\sup_{t\in[0,1]}\mathrm{var}(g_{n})_{t}\leq 1$, there exists a subsequence $(n_{k})_{k=1}^{\infty}$ and a finite variation $g$ such that $\lim_{k\to\infty}\int_{0}^{1}fdg_{n_{k}}=\int_{0}^{1}fdg$ for all continuous $f$. In Section \ref{sec:counter}, we show the continuity condition in Theorem \ref{thm:main} cannot be weakened—mirroring the deterministic case.
\par
The imposition of a boundedness condition on Theorem \ref{thm:main} is natural, and not a significant restriction. Indeed, the boundedness conditions assumed in the literature are often stronger (see, for example, \cite{filterappl,del-sch}). Furthermore, since the boundedness condition from Theorem \ref{thm:main} is based on boundedness in probability, Theorem \ref{thm:main} is not chained to any particular choice of measure, an important property when working with multiple equivalent probability measures.
\par
Let us note the following corollary of Theorem \ref{thm:main} for $H^{1}$-bounded sequences of martingales.
\begin{corollary}\label{corr:main-mart}
    Let $(M^{n})_{n=1}^{\infty}$ be a sequence of martingales such that
    \begin{equation*}
        \sup_{n}\int_{\Omega}\left(M^{n}\right)^{\ast}_{t}d\mathbb{P}<\infty
    \end{equation*}
    for each $t\geq0$. Then there exists $\widetilde{M}^{n}\in\mathrm{co}\{M^{m}:m\geq n\}$ and a semimartingale $\widetilde{M}$ such that $\widetilde{M}=\lim_{s\downarrow\cdot,s\in\mathbb{Q}_{+}}\lim_{n\to\infty}\widetilde{M}^{n}_{s}$ and
    \begin{equation*}
        \lim_{s\downarrow\cdot,s\in\mathbb{Q}_{+}}\mathbb{P}-\lim_{n\to\infty}\int_{0}^{s}Yd\widetilde{M}^{n}=\int_{0}^{\cdot}Yd\widetilde{M}
    \end{equation*}
    for each continuous semimartingale $Y$.
\end{corollary}
Note that one can say little about $\widetilde{M}$ in general; in particular, $\widetilde{M}$ may fail to be a martingale (see Example 3.1, \cite{del-sch}). However, if there exists an integrable random variable $\xi\geq0$ such that
    \begin{equation*}
        M^{n}_{t}\geq-\xi
    \end{equation*}
    for each $t\geq0$, then the Fatou limit $\widetilde{M}$ from Corollary \ref{corr:main-mart} is a supermartingale.
\begin{proof}
    Given the validity of Theorem \ref{thm:main}, it suffices to show that for each $t\geq0$, the set
    \begin{equation*}
        \mathrm{co}\left\{\left\vert{\int_{0}^{t}\xi dX^{n}}\right\vert:n\in\mathbb{N},\xi\textrm{ is predictable and }\vert{\xi}\vert\leq1\right\}
    \end{equation*}
    is bounded in probability. This is a simple consequence of the $p=1$ Burkholder-Davis-Gundy inequality.
\end{proof}
\section{Proof of main result}\label{sec:proof}
Our proof relies heavily on the following elementary consequence of integration by parts:
\begin{equation}\label{eq:ibp}
    \int_{0}^{\cdot}Y_{-}dZ=YZ-Y_{0}Z_{0}-\int_{0}^{\cdot}Z_{-}dY-[Z,Y]
\end{equation}
where $Z$ is any semimartingale (note that, since $Y$ is continuous, $Y=Y_{-}$). Previous attempts at compactness principles for sequences of the form $(\int_{0}^{\cdot}YdX^{n})_{n=1}^{\infty}$ rely to a great extent on such dual formulations of the problem (for example, see Proposition 2.12, \cite{martconvcomp}).
\par
Suppose we know \textit{a priori} that $\lim_{s\downarrow\cdot,s\in\mathbb{Q}_{+}}\lim_{n\to\infty}X^{n}_{s}$ exists $\mathbb{P}$-a.s., and is equal to a semimartingale $X$. When passing to a limit, the first two terms of $\int_{0}^{\cdot}YdX^{n}$ in the integration by parts formulation cause no problems. Any obstructions to Fatou convergence must therefore stem from the last two terms, which are difficult to tame. Indeed, when the continuity condition on $Y$ is dropped, we give a counterexample in \S\ref{sec:counter} showing that $([X^{n},Y])_{n=1}^{\infty}$ can fail to converge to the proper limit, even when all the other terms do. Another counterexample in \S\ref{sec:counter} for discontinuous $Y$ shows that the sum of all the terms in (\ref{eq:ibp}) can converge to the proper limit, even when none of the last two terms in (\ref{eq:ibp}) do. Our techniques therefore hinge on continuity of $Y$.
\par
There is no evidently clear global relationship between $X_{-}$ and $(X^{n}_{-})_{n=1}^{\infty}$, and hence also between $\int_{0}^{\cdot}X_{-}dY$ and $\left(\int_{0}^{\cdot}X^{n}_{-}dY\right)_{n=1}^{\infty}$. However, though $X_{-}$ may fail to be the pointwise limit of $(X^{n}_{-})_{n=1}^{\infty}$ (even when one merely requires convergence outside of an evanescent set), we show that, modulo many technicalities,\footnote{In particular, it may be necessary to replace the sequence $(X^{n})_{n=1}^{\infty}$ with a sequence of convex combinations thereof.} the latter does nevertheless approach the former on a predictable set $D$ whose complement is thin (in the sense of Definition 1.30, \cite{jacshir13}). Since $Y$ is continuous, such sets are ``large'' relative to $dY$. A stopping-time argument, controlling how large the relevant processes can get, allows one to conclude convergence of the integrals $\left(\int_{0}^{\cdot}X^{n}_{-}dY\right)_{n=1}^{\infty}$.
\par
For the quadratic covariation terms $([X^{n},Y])_{n=1}^{\infty}$, continuity of $Y$ is used to obtain that $[X,Y]$ and $([X^{n},Y])_{n=1}^{\infty}$ have no jumps. This is relevant, since for purely discontinuous $Y$ it is possible for $\left(\bigcup_{n=1}^{\infty}\{\Delta X^{n}\neq0\}\right)\cap\{\Delta Y\neq0\}$ to be evanescent, while
\begin{equation*}
    \mathbb{P}\left(\{\omega:\exists t\in[0,\infty)\textrm{ with }(\omega,t)\in\{\Delta X\neq0\}\cap\{\Delta Y\neq0\}\}\right)>0
\end{equation*}
so that $[X^{n},Y]=0$ for all $n$, but $[X,Y]\neq0$ (see Theorem \ref{thm:qlc} for an example). Furthermore, the assumption of continuity allows one to focus only on the noisy part of the relevant processes when making quadratic covariation calculations, which in questions of convergence are usually easier to deal with compared to compensator processes (see, for example, \cite{meminslum}).
\subsection{Lemmata}\label{subsec:lem}
By applying a localization procedure and diagonalization, we may assume the existence of a time $T\in[0,\infty)$ such that $X^{n}$ stopped at $T$ is $X^{n}$ for all $n$. Furthermore, it is no loss of generality to assume that $X^{n}_{0}=0$ for all $n$.
\begin{lemma}\label{lem:doobmeyer-switch}
    There exists an equivalent probability measure $\mathbb{Q}\sim\mathbb{P}$ and $\widetilde{X}^{n}\in\mathrm{co}\{X^{m}:m\geq n\}$ such that the following holds.
    \begin{enumerate}
        \item $\widetilde{X}^{n}$ is a $\mathbb{Q}$-special semimartingale for each $n$ with $\mathbb{Q}$-Doob-Meyer decomposition $\widetilde{X}^{n}=M^{n}+A^{n}$ (where $M^{n}$ is the local martingale part, and $A^{n}$ is the compensator).
        \item Each $M^{n}$ is an $L^{2}(\mathbb{Q})$-martingale, $(M^{n})_{n=1}^{\infty}$ converges in $L^{2}(\mathbb{Q})$ to some $L^{2}(\mathbb{Q})$-martingale $M$, and $\Vert{M^{n}-M}\Vert_{L^{2}(\mathbb{Q})}\leq\frac{1}{4^{n}}$.
        \item $\sup_{n}\int_{\Omega}\mathrm{var}(A^{n})_{\infty}d\mathbb{Q}<\infty$, and $\sup_{n}\mathrm{var}(A^{n})_{\infty}<\infty$ up to a $\mathbb{P}$-null set.
        \item $(A^{n})_{n=1}^{\infty}$ converges pointwise up to an evanescent set to some predictable finite variation process $A$.\footnote{By Froda's theorem (see Lemma \ref{lem:froda}), both the left limit process $A_{-}$ and right limit process $A_{+}$ are well-defined and adapted (for the former, this is obvious, while for the latter, it is a consequence of right-continuity of the filtration).}
        \item We have that $\widetilde{A}=\lim_{s\downarrow\cdot,s\in\mathbb{Q}_{+}}\lim_{n\to\infty}A^{n}_{s}$ exists $\mathbb{P}$-a.s. and is a finite variation semimartingale.
    \end{enumerate}
\end{lemma}
\begin{proof}
    We first apply some results by \cite{melnikovemery} for semimartingales on a finite time interval; they are applicable by our localization procedure (see the very beginning of \S\ref{subsec:lem}). By passing to convex combinations and applying (Lemma 7, \cite{melnikovemery}), we may assume that $\mathrm{co}\left\{[X^{n},X^{n}]_{\infty}:n\in\mathbb{N}\right\}$ is bounded in probability. Thus, using (Theorem 3, \cite{melnikovemery}), we can find a probability measure $\mathbb{Q}\sim\mathbb{P}$ such that each $X^{n}$ is $\mathbb{Q}$-special with Doob-Meyer decomposition $X^{n}=M^{n}+A^{n}$ for an $L^{2}(\mathbb{Q})$-martingale $M^{n}$ and a predictable finite variation process $A^{n}$, $(M^{n})_{n=1}^{\infty}$ forms a bounded sequence of $L^{2}(\mathbb{Q})$-martingales, and $(\mathrm{var}(A^{n})_{\infty})_{n=1}^{\infty}$ is bounded in $L^{1}(\mathbb{Q})$. Komlós's theorem implies, after a passage to convex combinations, that the latter part of (3) holds (see Lemma 2.5, \cite{schruessdies}).
    \par
    By applying (Theorem 3.19, \cite{brezis}) in tandem with (Corollary 3.8, \cite{brezis}) to the Hilbert space of $L^{2}(\mathbb{Q})$-martingales and passing to convex combinations, we may assume that $(M^{n})_{n=1}^{\infty}$ converges in $L^{2}(\mathbb{Q})$ to some $M$, and $\Vert{M^{n}-M}\Vert_{L^{2}(\mathbb{Q})}\leq\frac{1}{4^{n}}$ (showing (2)).
    \par
    By (Proposition 3.4, \cite{schcamp}) we have (4) after passing to convex combinations. Passing to further convex combinations and applying (Lemma 2.7, \cite{del-sch}) and (Remark 2.8, \cite{del-sch}), we obtain (5).
\end{proof}
\begin{lemma}\label{lem:somethinsets}
    The sets $\{A\neq A_{-}\}$, $\{A\neq A_{+}\}$, and $\{A_{-}\neq\widetilde{A}_{-}\}$ are thin.
\end{lemma}
\begin{proof}
    By Froda's theorem (see Lemma \ref{lem:froda}), the optional sets $\{A\neq A_{-}\}$ and $\{A\neq A_{+}\}$ have at most countably many sections, and are therefore thin by Theorem 117 in Appendix IV of \cite{meyerdell}. Remark also that $\{\widetilde{A}\neq\widetilde{A}_{-}\}$ is thin (use a similar argument to the one we just gave, or see Proposition 1.32, \cite{jacshir13}), a fact we will need for the next paragraph.
    \par
    We will now prove thinness of $\{A_{-}\neq\widetilde{A}_{-}\}$. It suffices to show that, modulo an evanescent set, $\{A_{-}\neq\widetilde{A}_{-}\}$ is contained in a thin set. Let
    \begin{equation*}
        E=\{A\neq A_{-}\}\cup\{A\neq A_{+}\}\cup\{\widetilde{A}\neq\widetilde{A}_{-}\}
    \end{equation*}
    which is thin by the previous paragraph. It is clear that $\{A_{-}\neq\widetilde{A}_{-}\}$ is contained in $E$ (modulo an evanescent set), which proves the claim.
\end{proof}
\begin{lemma}\label{lem:integralonthin}
    Let $Z$ be a quasi-left-continuous semimartingale, and let $G$ and $H$ be $Z$-integrable predictable processes. If $\{G\neq H\}$ is thin, then
    \begin{equation*}
        \int_{0}^{\cdot}GdZ=\int_{0}^{\cdot}HdZ
    \end{equation*}
    up to evanescence.
\end{lemma}
\begin{proof}
    Since $G$ and $H$ are predictable, $\{G\neq H\}$ is a predictable set. Thus, by thinness and (Theorem IV.81.c, \cite{meyerdell}), we may write
    \begin{equation*}
        \{G\neq H\}=\bigcup_{n=1}^{\infty}\llbracket{\sigma_{n}}\rrbracket,
    \end{equation*}
    for some sequence $(\sigma_{n})_{n=1}^{\infty}$ of predictable stopping times with pairwise disjoint graphs. Let $E_{n}=\bigcup_{m=n}^{\infty}\llbracket{\sigma_{m}}\rrbracket$. Then
    \begin{equation*}
        \int_{0}^{\cdot}\mathbf{1}_{(\Omega\times[0,\infty))\setminus E_{n}}(G-H)dZ=\int_{0}^{\cdot}\mathbf{1}_{\{G=H\}}(G-H)dZ+\sum_{i=1}^{n-1}(G_{\sigma_{i}}-H_{\sigma_{i}})\Delta Z_{\sigma_{i}}\mathbf{1}_{\llbracket{\sigma_{i},\infty}\llbracket}=0
    \end{equation*}
    since $Z$ is quasi-left-continuous, and $\sigma_{i}$ is predictable. Taking $n\to\infty$, we therefore obtain that
    \begin{equation*}
        0=\lim_{n\to\infty}\int_{0}^{\cdot}\mathbf{1}_{(\Omega\times[0,\infty))\setminus E_{n}}(G-H)dZ=\int_{0}^{\cdot}(G-H)dZ
    \end{equation*}
    in the Émery topology by the bounded convergence theorem, so that $\int_{0}^{\cdot}(G-H)dZ=0$, as desired.
\end{proof}
\subsection{The proof of Theorem \ref{thm:main}}
\begin{proof}[Proof of Theorem \ref{thm:main}]
    Let $\widetilde{X}=M+\widetilde{A}$. We first show that $\lim_{s\downarrow\cdot}\lim_{n\to\infty}\widetilde{X}^{n}_{s}$ exists $\mathbb{P}$-a.s. and is equal to $\widetilde{X}$. It suffices to show that $((M^{n}-M)^{\ast}_{\infty})_{n=1}^{\infty}$ $\mathbb{P}$-a.s. converges to zero. Notice that, by Markov's inequality and Doob's maximal inequality,
    \begin{equation*}
        \mathbb{Q}\left(\left\{(M^{n}-M)^{\ast}_{\infty}\geq\frac{1}{2^{n}}\right\}\right)\leq2^{n}\int_{\Omega}(M^{n}-M)^{\ast}_{\infty}d\mathbb{Q}
    \end{equation*}
    \begin{equation*}
        \leq2^{n+1}\Vert{M^{n}-M}\Vert_{L^{2}(\mathbb{Q})}\leq2^{n+1}\frac{1}{4^{n}}=2^{1-n}.
    \end{equation*}
    The Borel-Cantelli lemma implies that $\mathbb{Q}\left(\limsup_{n\to\infty}\left\{(M^{n}-M)^{\ast}_{\infty}\geq\frac{1}{2^{n}}\right\}\right)=0$; the equivalence of $\mathbb{Q}$ and $\mathbb{P}$ shows that $\mathbb{P}\left(\limsup_{n\to\infty}\left\{(M^{n}-M)^{\ast}_{\infty}\geq\frac{1}{2^{n}}\right\}\right)=0$, implying $\mathbb{P}$-a.s. convergence of $((M^{n}-M)^{\ast}_{\infty})_{n=1}^{\infty}$ to zero.
    \par
    By equation (\ref{eq:ibp}), for the validity of Theorem \ref{thm:main}, it is enough to show that
    \begin{equation*}
        \lim_{s\downarrow\cdot}\mathbb{P}-\lim_{n\to\infty}\left(Y_{s}\widetilde{X}^{n}_{s}-\int_{0}^{s}\widetilde{X}^{n}_{-}dY-[\widetilde{X}^{n},Y]_{s}\right)=Y\widetilde{X}-\int_{0}^{\cdot}\widetilde{X}_{-}dY-[\widetilde{X},Y].
    \end{equation*}
    Thus, it suffices to show (1)$\wedge$(2)$\wedge$(3), where:
    \begin{enumerate}
        \item $\lim_{s\downarrow\cdot}\lim_{n\to\infty}Y_{s}\widetilde{X}^{n}_{s}=Y\widetilde{X}$ exists $\mathbb{P}$-a.s. and equals $Y\widetilde{X}$.
        \item $\left(\int_{0}^{\cdot}\widetilde{X}^{n}_{-}dY-\int_{0}^{\cdot}\widetilde{X}_{-}dY\right)_{n=1}^{\infty}$ converges to zero in the Émery topology (in particular, in the u.c.p. topology).
        \item  $\left(\mathrm{var}\left([\widetilde{X}^{n},Y]-[\widetilde{X},Y]\right)_{t}\right)_{n=1}^{\infty}$ converges $\mathbb{P}$-a.s. to zero for each $t\geq0$.
    \end{enumerate}
    (1) is clear from $\widetilde{X}=\lim_{s\downarrow\cdot}\lim_{n\to\infty}\widetilde{X}^{n}_{s}$.
    \par
    For (2), by splitting $\widetilde{X}^{n}_{-}$ as $\widetilde{X}^{n}_{-}=M^{n}_{-}+A^{n}_{-}$ and $\widetilde{X}$ as $\widetilde{X}=M+\widetilde{A}$, it suffices to show that $\left(\int_{0}^{\cdot}M^{n}_{-}dY-\int_{0}^{\cdot}M_{-}dY\right)_{n=1}^{\infty}$ and $\left(\int_{0}^{\cdot}A^{n}_{-}dY-\int_{0}^{\cdot}\widetilde{A}_{-}dY\right)_{n=1}^{\infty}$ converge to zero in the Émery topology. 
    \par
    We first show Émery convergence of $\left(\int_{0}^{\cdot}M^{n}_{-}dY-\int_{0}^{\cdot}M_{-}dY\right)_{n=1}^{\infty}$ to zero. Let $G=\sup_{n}(M^{n}-M)^{\ast}$ (which is finite up to an evanescent set, by $\mathbb{P}$-a.s. convergence of $((M^{n}-M)^{\ast}_{\infty})_{n=1}^{\infty}$), and $H=G_{-}$. It is not difficult to see that $H$ is $Y$-integrable. By $\mathbb{P}$-a.s. convergence of $((M^{n}-M)^{\ast}_{\infty})_{n=1}^{\infty}$ to zero, $(M^{n}_{-}-M_{-})_{n=1}^{\infty}$ converges pointwise (modulo evanescent sets) to zero. Since $\vert{M^{n}_{-}-M_{-}}\vert\leq\vert{H}\vert$ for each $n$, the stochastic dominated convergence theorem therefore shows that $\left(\int_{0}^{\cdot}M^{n}_{-}dY-\int_{0}^{\cdot}M_{-}dY\right)_{n=1}^{\infty}$ Émery-converges to zero.
    \par
    We now show Émery convergence of $\left(\int_{0}^{\cdot}A^{n}_{-}dY-\int_{0}^{\cdot}\widetilde{A}_{-}dY\right)_{n=1}^{\infty}$ to zero. $A$ is $Y$-integrable, and $A^{n}$ is $Y$-integrable for each $n$. Let $E=\bigcup_{n=1}^{\infty}\{A^{n}\neq A^{n}_{-}\}$, which is predictable and thin, as $\{A^{n}\neq A^{n}_{-}\}$ is predictable and thin for each $n$ (see Proposition 1.32, \cite{jacshir13}). Let $J=\sup_{n}\mathrm{var}(A^{n})$, and define a sequence $(\sigma_{m})_{m=1}^{\infty}$ of stopping times by $\sigma_{m}=\inf\{t:J_{t}\geq m\}$. Since $\sup_{n}\mathrm{var}(A^{n})_{\infty}<\infty$ $\mathbb{P}$-a.s., $(\sigma_{m})_{m=1}^{\infty}$ is a localizing sequence. We have that $\left\vert\left(\mathbf{1}_{(\Omega\times[0,\infty))\setminus E}A^{n}\right)^{\sigma_{m}}\right\vert\leq m$ for each $m$ and $n$. Thus, since $(A^{n})_{n=1}^{\infty}$ converges to $A$ pointwise, the stochastic dominated convergence theorem implies that $\left(\int_{0}^{\cdot}\mathbf{1}_{(\Omega\times[0,\infty))\setminus E}A^{n}dY\right)_{n=1}^{\infty}$ converges to $\int_{0}^{\cdot}\mathbf{1}_{(\Omega\times[0,\infty))\setminus E}AdY$ in the Émery topology. By Lemma \ref{lem:integralonthin} and Lemma \ref{lem:somethinsets}, we have that 
    \begin{equation*}
        \forall n:\int_{0}^{\cdot}\mathbf{1}_{(\Omega\times[0,\infty))\setminus E}A^{n}dY=\int_{0}^{\cdot}A^{n}dY=\int_{0}^{\cdot}A^{n}_{-}dY,
    \end{equation*}
    \begin{equation*}\int_{0}^{\cdot}\mathbf{1}_{(\Omega\times[0,\infty))\setminus E}AdY=\int_{0}^{\cdot}AdY=\int_{0}^{\cdot}A_{-}dY=\int_{0}^{\cdot}\widetilde{A}_{-}dY.
    \end{equation*} 
    Thus, $\left(\int_{0}^{\cdot}A^{n}_{-}dY\right)_{n=1}^{\infty}$ converges in the Émery topology to $\int_{0}^{\cdot}\widetilde{A}_{-}dY$, as desired. This finishes the proof of (2).
    \par
    We now show (3). $Y$ is a special $\mathbb{Q}$-semimartingale, as it is continuous; write the $\mathbb{Q}$-Doob-Meyer decomposition of $Y$ as $Y=Y_{0}+N+B$, where $N$ is a continuous local martingale with $N_{0}=0$, and $B$ is a continuous finite variation semimartingale. For any semimartingale $Z$ which decomposes as $Z=L+V$, where $L$ is a local martingale, and $V$ is a finite variation semimartingale, we have
    \begin{equation*}\label{eq:quadcancelmart}
        [Z,Y]=[L,N]+\sum_{s\leq\cdot}\left(\Delta L_{s}\Delta B_{s}+\Delta V_{s}\Delta N_{s}+\Delta V_{s}\Delta B_{s}\right)=[L,N].
    \end{equation*}
    Thus, it suffices to show that $\left(\mathrm{var}\left([M^{n},N]-[M,N]\right)_{t}\right)_{n=1}^{\infty}$ converges $\mathbb{P}$-a.s. to zero for each $t\geq0$.
    \par
    Let $\tau_{n}=\inf\{t:\vert{N_{t}}\vert\geq n\}$; it is clear that $(\tau_{n})_{n=1}^{\infty}$ is a localizing sequence. It suffices to show that $\left(\mathrm{var}\left([M^{n},N]-[M,N]\right)_{\tau_{m}}\right)_{n=1}^{\infty}$ converges $\mathbb{P}$-a.s. to zero for each $m$. By Markov's inequality
    \begin{equation*}
        \mathbb{Q}\left(\left\{\mathrm{var}\left([M^{n},N]-[M,N]\right)_{\tau_{m}}\geq\frac{1}{2^{n}}\right\}\right)\leq2^{n}\int_{\Omega}\mathrm{var}\left([M^{n}-M,N]\right)_{\tau_{m}}d\mathbb{Q}.
    \end{equation*}
    By the Kunita-Watanabe inequality,
    \begin{equation*}
        \int_{\Omega}\mathrm{var}\left([M^{n}-M,N]\right)_{\tau_{m}}d\mathbb{Q}\leq\Vert{M^{n}-M}\Vert_{L^{2}(\mathbb{Q})}\Vert{N^{\tau_{m}}}\Vert_{L^{2}(\mathbb{Q})}\leq\frac{m}{4^{n}}.
    \end{equation*}
    Thus,
    \begin{equation*}
       \mathbb{Q}\left(\left\{\mathrm{var}\left([M^{n},N]-[M,N]\right)_{\tau_{m}}\geq\frac{1}{2^{n}}\right\}\right)\leq\frac{m2^{n}}{4^{n}}=\frac{m}{2^{n}}.
    \end{equation*}
    The Borel-Cantelli lemma implies that $\limsup_{n\to\infty}\left\{\mathrm{var}\left([M^{n},N]-[M,N]\right)_{\tau_{m}}\geq\frac{1}{2^{n}}\right\}$ has $\mathbb{Q}$-measure zero, hence also $\mathbb{P}$-measure zero. Thus, $\left(\mathrm{var}\left([M^{n},N]-[M,N]\right)_{\tau_{m}}\right)_{n=1}^{\infty}$ converges $\mathbb{P}$-a.s. to zero for each $m$, as desired. This completes the proof of Theorem \ref{thm:main}.
\end{proof}
\begin{remark}
    For a fixed $Y$, it is of course possible via the Borel-Cantelli lemma to pass to a subsequence $(n_{k})_{k=1}^{\infty}$ such that
    \begin{equation*}
        \lim_{s\downarrow\cdot,s\in\mathbb{Q}_{+}}\lim_{k\to\infty}\int_{0}^{s}Yd\widetilde{X}^{n_{k}}=\int_{0}^{\cdot}Yd\widetilde{X}.
    \end{equation*}
    It is unclear whether one can choose $(n_{k})_{k=1}^{\infty}$ so that the above holds for all $Y$ simultaneously.
\end{remark}
\section{Necessity of a continuous integrand}\label{sec:counter}
Theorem \ref{thm:main} assumes an integrand which is continuous and a semimartingale. The integration by parts techniques used to prove Theorem \ref{thm:main} hint at generalizations of Theorem \ref{thm:main} for integrands of the form $Y_{-}$ for a general semimartingale $Y$ (especially given that a semimartingale $Y$ need not be an admissible integrand in the Itô theory, if $Y$ is not continuous). 
\par
We show that Theorem \ref{thm:main} is sharp, in the sense that one cannot hope to extend Theorem \ref{thm:main} beyond continuous semimartingale integrands. More precisely, we have the following.
\begin{theorem}\label{thm:qlc}
    There exists a semimartingale $Y$, and a sequence $(X^{n})_{n=1}^{\infty}$ of semimartingales starting at zero such that the following holds.
    \begin{enumerate}
        \item $\sup_{n}\int_{\Omega}\mathrm{var}(X^{n})_{\infty}d\mathbb{P}\leq1$.
        \item If $\widetilde{X}^{n}\in\mathrm{co}\{X^{m}:m\geq n\}$, then $Y=\lim_{s\downarrow\cdot,s\in\mathbb{Q}_{+}}\lim_{n\to\infty}\widetilde{X}^{n}_{s}$.
        \item If $\widetilde{X}^{n}\in\mathrm{co}\{X^{m}:m\geq n\}$, then
        \begin{equation*}
        \lim_{s\downarrow\cdot,s\in\mathbb{Q}_{+}}\lim_{n\to\infty}\int_{0}^{s}Y_{-}d\widetilde{X}^{n}_{-}\neq\int_{0}^{\cdot}Y_{-}dY,
        \end{equation*}
    \end{enumerate}
\end{theorem}
\begin{proof}
    Define $Y=\mathbf{1}_{\llbracket{1,\infty}\llbracket}$, and let $X^{n}=\mathbf{1}_{\llbracket{1+1/n,\infty}\llbracket}$. If $\widetilde{X}^{n}\in\mathrm{co}\{X^{m}:m\geq n\}$, we have that $Y=\lim_{s\downarrow\cdot,s\in\mathbb{Q}_{+}}\lim_{n\to\infty}\widetilde{X}^{n}_{s}$ and $\lim_{s\downarrow\cdot,s\in\mathbb{Q}_{+}}\lim_{n\to\infty}\int_{0}^{s}Y_{-}d\widetilde{X}^{n}_{-}=\mathbf{1}_{\llbracket{1,\infty}\llbracket}\neq0$, but $\int_{0}^{\cdot}Y_{-}dY=0$.
\end{proof}
\begin{remark}
    If $\mathbb{F}$ admits a totally inaccessible stopping time, the above counterexample can be modified so that $Y$ is quasi-left-continuous (c.f. Remark 2.9, \cite{del-sch}).
\end{remark}
However, even when $\left(\int_{0}^{\cdot}Y_{-}dX^{n}\right)_{n=1}^{\infty}$ converges to the proper limit, it is possible that the techniques used in \S\ref{sec:proof}—the integration by parts formula (\ref{eq:ibp})—are of no use.
\begin{theorem}\label{thm:convg-but-counter}
    There exists a semimartingale $Y$, and a sequence $(X^{n})_{n=1}^{\infty}$ of semimartingales starting at zero such that the following holds.
    \begin{enumerate}
        \item $\sup_{n}\int_{\Omega}\mathrm{var}(X^{n})_{\infty}d\mathbb{P}\leq1$.
        \item If $\widetilde{X}^{n}\in\mathrm{co}\{X^{m}:m\geq n\}$, then $Y=\lim_{s\downarrow\cdot,s\in\mathbb{Q}_{+}}\lim_{n\to\infty}\widetilde{X}^{n}_{s}$.
        \item If $\widetilde{X}^{n}\in\mathrm{co}\{X^{m}:m\geq n\}$, then
        \begin{equation*}
            \lim_{s\downarrow\cdot,s\in\mathbb{Q}_{+}}\lim_{n\to\infty}\int_{0}^{s}Y_{-}d\widetilde{X}^{n}_{-}=\int_{0}^{\cdot}Y_{-}dY,
        \end{equation*}
        \begin{equation*}
            \lim_{s\downarrow\cdot,s\in\mathbb{Q}_{+}}\lim_{n\to\infty}[\widetilde{X}^{n},Y]_{s}\neq[Y,Y],
        \end{equation*}
        \begin{equation*}
            \lim_{s\downarrow\cdot,s\in\mathbb{Q}_{+}}\lim_{n\to\infty}\int_{0}^{s}\widetilde{X}^{n}_{-}dY\neq\int_{0}^{\cdot}Y_{-}dY.
        \end{equation*}
    \end{enumerate}
\end{theorem}
\begin{proof}
    Define $Y=\mathbf{1}_{\llbracket{1,\infty}\llbracket}$, and let $X^{n}=\mathbf{1}_{\llbracket{1-1/n,\infty}\llbracket}$ (c.f. Example 1, \cite{jakcounter}). If $\widetilde{X}^{n}\in\mathrm{co}\{X^{m}:m\geq n\}$, we have that $Y=\lim_{s\downarrow\cdot,s\in\mathbb{Q}_{+}}\lim_{n\to\infty}\widetilde{X}^{n}_{s}$ and $\int_{0}^{\cdot}Y_{-}d\widetilde{X}^{n}_{-}=0=\int_{0}^{\cdot}Y_{-}dY$. However, $[Y,Y]\neq0=[\widetilde{X}^{n},Y]$ and $\lim_{s\downarrow\cdot,s\in\mathbb{Q}_{+}}\lim_{n\to\infty}\int_{0}^{s}\widetilde{X}^{n}_{-}dY\neq0=\int_{0}^{\cdot}Y_{-}dY$.
\end{proof}
\printbibliography
\end{document}